\documentclass[12pt]{article}
\usepackage[centertags]{amsmath}
\usepackage{amsfonts}
\usepackage{amssymb}
\usepackage{amsthm}
\usepackage{newlfont}

\numberwithin{equation}{section}

\newtheorem{theorem}{Theorem}

\newtheorem{definition}{Definition}

\begin{document}

\title{\textbf{ On the non-Abelian Lie bracket and the generalized covariant Hamilton system}}
\vspace{.3 cm}
\author{Gen Wang\footnote{ \small Department of Mathematics, Zhejiang Normal University, Zhejiang, Jinhua 321004, China. Email: wanggen@zjnu.edu.cn}}
\date{}

\maketitle

\begin{abstract}
Based on the non-Abelian Lie algebra, a generalized geometric Lie bracket on vector space is proposed to further realize the generalized structural Poisson bracket, and then we briefly discuss the second order equations of the generalized covariant Hamilton system with respect to the time and coordinates and its simple application.

\end{abstract}


\section{Preliminaries}
\subsection{Non-Abelian Lie algebra}

A non-abelian group is a group $\left( \mathcal{G},* \right)$ in which there exists at least one pair of elements $a$ and $b$ of $\mathcal{G}$, such that $ a*  b \neq b *  a$.  Both discrete groups and continuous groups may be non-abelian.  A Lie algebra is a vector space ${\mathfrak {g}}$ together with a non-associative, alternating bilinear map $${\displaystyle {\mathfrak {g}}\times {\mathfrak {g}}\rightarrow {\mathfrak {g}}}$$  $${(\mathbf{x},\mathbf{y})\mapsto [\mathbf{x},\mathbf{y}]_{cr}}=\mathbf{x}* \mathbf{y}-\mathbf{y}* \mathbf{x}\in {\mathfrak {g}},$$ called the Lie bracket, satisfying the Jacobi identity.   For any associative algebra $L$ with multiplication $*$, The Lie bracket of two elements of Lie algebra is defined to be their commutator in $L$:  $[a,b]_{cr}=a*b-b*a$.  The associativity of the multiplication $*$ in $L$ implies the Jacobi identity of the commutator in Lie bracket.

\begin{definition}\cite{1,2}\label{de1}
A Lie algebra is a vector space ${\mathfrak {g}}$ over some field $\mathbb{F}$ together with a binary operation $ [\cdot ,\cdot ]_{cr}:{\mathfrak {g}}\times {\mathfrak {g}}\to {\mathfrak {g}}$ called the Lie bracket that satisfies the following axioms:
\begin{description}
  \item[(i)] Bilinearity,
$[a_{1}\mathbf{x}+a_{2}\mathbf{y},\mathbf{z}]_{cr}=a_{1}[\mathbf{x},\mathbf{z}]_{cr}+a_{2}[\mathbf{y},\mathbf{z}]_{cr},\quad [\mathbf{z},a_{1}\mathbf{x}+a_{2}\mathbf{y}]_{cr}=a_{1}[\mathbf{z},\mathbf{x}]_{cr}+a_{2}[\mathbf{z},\mathbf{y}]_{cr}$
for all scalars $a_{1}, a_{2}\in \mathbb{F}$ and all elements $\mathbf{x},\mathbf{y}, \mathbf{z}\in {\mathfrak {g}}$.
  \item[(ii)]
Alternativity,
$[\mathbf{x},\mathbf{x}]_{cr}=0$, for all $\mathbf{x}\in {\mathfrak {g}}$.
  \item[(iii)]
The Jacobi identity,
$[\mathbf{x},[\mathbf{y},\mathbf{z}]_{cr}]_{cr}+[\mathbf{z},[\mathbf{x},\mathbf{y}]_{cr}]_{cr}+[\mathbf{y},[\mathbf{z},\mathbf{x}]_{cr}]_{cr}=0$,
for all $\mathbf{x},\mathbf{y}, \mathbf{z}\in {\mathfrak {g}}$.
  \item[(iv)]
Anticommutativity,  ${\displaystyle [\mathbf{x},\mathbf{y}]_{cr}=-[\mathbf{y},\mathbf{x}]_{cr} }$,
for all elements $\mathbf{x},\mathbf{y}\in {\mathfrak {g}}$.

\end{description}
\end{definition}
It is habitual to express a Lie algebra like ${\mathfrak {g}}$, and if a Lie algebra is associated with a Lie group, then the spelling of the Lie algebra is the same as that Lie group.  Most of the interesting Lie groups are non-abelian, 
\begin{theorem}\cite{1}
  Let ${{X}_{1}},\cdots ,{{X}_{n}}$ be a basis of Lie algebra $\mathcal{G}$ of Lie groups, $\left[ {{X}_{i}},{{X}_{j}} \right]_{cr}\in \mathcal{G}$, then there exists structural constants $c_{ij}^{k}\in {{C}^{\infty }}\left( M,\mathbb{R} \right)$ such that $\left[ {{X}_{i}},{{X}_{j}} \right]_{cr}=c_{ij}^{k}{{X}_{k}}$ satisfying the following properties
\begin{enumerate}
  \item $c_{ij}^{k}+c_{ji}^{k}=0$.
  \item $c_{ij}^{r}c_{rk}^{s}+c_{jk}^{r}c_{ri}^{s}+c_{ki}^{r}c_{rj}^{s}=0$.
\end{enumerate}
\end{theorem}
Suppose that $\frac{\partial }{\partial {{x}^{1}}},\cdots ,\frac{\partial }{\partial {{x}^{n}}}$  is the basis of  Lie algebra $\mathcal{G}$ on the Euclidean space ${{\mathbb{R}}^{n}}$, then it yields $\left[ {{X}_{i}},{{X}_{j}} \right]_{cr}=0$, where ${{X}_{i}}=\frac{\partial }{\partial {{x}^{i}}}\in \mathcal{G}$, subsequently, structure constant $c_{ij}^{k}\in {{C}^{\infty }}\left( M,\mathbb{R} \right)$  of Lie groups hold $c_{ij}^{k}=0$.  

\subsection{Generalized geometric Lie bracket}
As the Lie algebra defined above, we can consider a element $\xi\in \mathfrak {g}$ that has nothing to do with the elements $\mathbf{x},\mathbf{y}\in \mathfrak {g}$, where $\xi\in \mathfrak {g}$ is completely determined by the structural property of the vector space $\mathfrak {g}$, this feature belongs to the intrinsic property of  $\mathfrak {g}$ itself. Therefore, as the generalized structural Poisson bracket defined, we can define a geometric bracket on the vector space $\mathfrak {g}$,
$$G\left( \xi ,\mathbf{x},\mathbf{y}\right)=\mathbf{x}{{\left[\xi,\mathbf{y} \right]}_{cr}}-\mathbf{y}{{\left[ \xi ,\mathbf{x}\right]}_{cr}},~~~\mathbf{x},\mathbf{y},\xi\in \mathfrak {g}$$
It can be seen that the geometric bracket on the vector space $\mathfrak {g}$ satisfies the anti-symmetry $G\left( \xi ,\mathbf{x},\mathbf{y} \right)=-G\left(\xi ,\mathbf{y},\mathbf{x}\right)$, if $\xi=0$ is given, then geometric bracket on the vector space $\mathfrak {g}$ disappears, that is $G\left( 0 ,\mathbf{x},\mathbf{y}\right)=0$.

\begin{definition}
A generalized geometric Lie algebra is a vector space ${\mathfrak {g}}$ over some field $\mathbb{F}$  together with a binary operation $ [\cdot ,\cdot ]=[\cdot ,\cdot ]_{cr}+G(\xi,\cdot ,\cdot ):{\mathfrak {g}}\times {\mathfrak {g}}\to {\mathfrak {g}}$ called the generalized geometric Lie bracket that satisfies the following axioms:
\begin{description}
  \item[(i)] Bilinearity,
$[a_{1}\mathbf{x}+a_{2}\mathbf{y},\mathbf{z}]=a_{1}[\mathbf{x},\mathbf{z}]+a_{2}[\mathbf{y},\mathbf{z}],\quad [\mathbf{z},a_{1}\mathbf{x}+a_{2}\mathbf{y}]=a_{1}[\mathbf{z},\mathbf{x}]+a_{2}[\mathbf{z},\mathbf{y}]$
for all scalars $a_{1}, a_{2}\in \mathbb{F}$ and all elements $\mathbf{x},\mathbf{y}, \mathbf{z}\in {\mathfrak {g}}$.
  \item[(ii)]
Alternativity,
$[\mathbf{x},\mathbf{x}]=0$, for all $\mathbf{x}\in {\mathfrak {g}}$.
  \item[(iii)]
The generalized Jacobi identity,
$[\mathbf{x},[\mathbf{y},\mathbf{z}]]+[\mathbf{z},[\mathbf{x},\mathbf{y}]]+[\mathbf{y},[\mathbf{z},\mathbf{x}]]=0$,
for all $\mathbf{x},\mathbf{y}, \mathbf{z}\in {\mathfrak {g}}$.
  \item[(iv)]
Anticommutativity,  ${\displaystyle [\mathbf{x},\mathbf{y}]=-[\mathbf{y},\mathbf{x}] }$,
for all elements $\mathbf{x},\mathbf{y}\in {\mathfrak {g}}$.

\end{description}

\end{definition}
Hence, the generalized geometric Lie bracket is completely written as 
\[\left[ \mathbf{x},\mathbf{y} \right]={{\left[ \mathbf{x},\mathbf{y} \right]}_{cr}}+G\left( \xi ,\mathbf{x},\mathbf{y}\right),~~\xi ,\mathbf{x},\mathbf{y}\in \mathfrak {g}\]
Therefore, if $\frac{\partial }{\partial {{x}^{1}}},\cdots ,\frac{\partial }{\partial {{x}^{n}}}$  is the basis of Lie algebra $\mathcal{G}$ on the Euclidean space ${{\mathbb{R}}^{n}}$, it can always choose a proper $\xi={{\xi}^{j}}{{\partial }_{j}}\in \mathcal{G}$ such that $G\left( \xi ,\mathbf{x},\mathbf{y}\right)\in \mathcal{G}$.

\section{Generalized structural Poisson bracket}
In this section, we briefly retrospect some basic concepts of different Poisson bracket, specifically \cite{3,4,5,6,7}. 

The classical Poisson bracket (CPB) defined on functions on ${{\mathbb{R}}^{2n}}$ is
\[\left\{ {{f}_{1}},{{f}_{2}} \right\}_{CPB}= \sum\limits_{i}{\left(\frac{\partial {{f}_{1}}}{\partial {{q}_{i}}}\frac{\partial {{f}_{2}}}{\partial {{p}_{i}}}-\frac{\partial {{f}_{1}}}{\partial {{p}_{i}}}\frac{\partial {{f}_{2}}}{\partial {{q}_{i}}} \right)},~~~\forall {{f}_{j}}\in {{C}^{\infty }}\left( M,\mathbb{R} \right)\]
By using classical Poisson bracket, it follows the classical Hamilton system
$$
  {\dot {q}}_{i}=\left\{ {{q}_{i}},H \right\}_{CPB},~~{\dot {p}}_{i}=\left\{ {{p}_{i}},H \right\}_{CPB},~~i=1,\cdots ,n$$
where ${\dot {q}}_{i}=\frac{d}{dt}{{q}_{i}}$. The variable $t$
denotes a real scalar variable called time, and the symbol
$\cdot$ is used for $d/dt$.  Meanwhile, the
canonical commutation relation in terms of the classical Poisson bracket follows
$\left\{ {q}_{j},{p}_{i} \right\}_{CPB}= \delta_{ij}$.

The generalized Poisson bracket (GPB) on ${{\mathbb{R}}^{r}}$ that is defined as the bilinear operation
\[{{\left\{ f,g \right\}}_{GPB}}={{\nabla }^{T}}fJ\nabla g={{J}_{ij}}\frac{\partial f}{\partial {{x}_{i}}}\frac{\partial g}{\partial {{x}_{j}}},~~i,j=1,2,\cdots n\] where structural matrix $J$ satisfies antisymmetric ${{J}_{ij}}={{\left\{ {{x}_{i}},{{x}_{j}} \right\}}_{GPB}}=-{{J}_{ji}}$. The generalized Hamiltonian system (GHS) is given by
\[{\dot{x}}=\frac{dx}{dt}=J\left( x \right)\nabla H\left( x \right),~~~x\in {{\mathbb{R}}^{n}}\]
where $\nabla H\left( x \right)$  is the gradient of the function Hamilton. Equivalently, using the generalized Poisson bracket above, the Hamiltonian equation can be further written as $${\dot{x}_{i}}=\left\{ {{x}_{i}},H \right\}_{GPB}={{{J}_{ij}}\frac{\partial H}{\partial {{x}_{j}}}}.$$The generalized
Hamiltonian system has widely been applied in many fields, but it still faces some difficulties. 

Later, we try to extend it to work for more general cases,
the generalized structural Poisson bracket (GSPB) on ${{\mathbb{R}}^{r}}$ of two functions $f,g\in {{C}^{\infty }}\left( M,\mathbb{R} \right)$ is shown as
\[\left\{ f,g \right\}={{\left\{ f,g \right\}}_{GPB}}+G\left(s, f,g \right), \]where $s\in {{C}^{\infty }}\left( M,\mathbb{R} \right)$ is the geometric potential function only given by the manifold space itself,  sometimes, it can be treated as the external environment, and it's generally chosen a form of logarithmic function,
and geometric bracket is $$G\left(s,f,g \right) =f{{\left\{ s ,g \right\}}_{GPB}}-g{{\left\{ s ,f \right\}}_{GPB}}.$$It implies the interaction with the external environment, 
where $\left\{ f,g \right\}=-\left\{ g,f \right\}$ is skew-symmetric.

By using the generalized gradient operator $D=\nabla +A$, and $A=\nabla s$,  the generalized structural Poisson bracket in terms of two functions $f,g\in {{C}^{\infty }}\left( M,\mathbb{R} \right)$ is also given by an analytic expression $$\left\{ f,g \right\}={{D}^{T}}fJDg={{J}_{ij}}{{D}_{i}}f{{D}_{j}}g\in {{C}^{\infty }}\left( M,\mathbb{R} \right),$$
where ${D}_{i}=\partial_{i}+{A}_{i}$, and ${A}_{i}=\partial_{i}s$ is structural derivative.

\begin{theorem}\label{th3}
For all $f,g,h \in  {{C}^{\infty }}\left( M,\mathbb{R} \right)$, $\lambda,\mu \in \mathbb{R}$, the generalized structural Poisson bracket has the following important properties

(1) Symmetry: $\left\{ f,g \right\}=-\left\{ g,f \right\}$.

(2) Bilinearity: $\left\{ \lambda f+\mu g,h \right\}=\lambda \left\{ f,h \right\}+\mu \left\{ g,h \right\}$.

(3) Generalized Jacobi identity: $\left\{ f,\left\{ g,h \right\} \right\}+\left\{ g,\left\{ h,f \right\} \right\}+\left\{ h,\left\{ f,g \right\} \right\}=0$.

(4) Generalized Leibnitz identity:  ${{\left\{ fg,h \right\}}}={{\left\{ fg,h \right\}}_{GPB}}+G\left(s, fg, h \right) $.

(5) Non degeneracy: if for all $f$, $\left\{ f,g \right\}=0$, then ${{\left\{ f,g \right\}}_{GPB}}=G\left(s, g,f \right)$.

\end{theorem}
Note that for the non degeneracy, $\left\{ f,g \right\}=0$ holds for all $f$,  
then it has an equation ${{\left\{ f,g \right\}}_{GPB}}=g{{\left\{ s,f \right\}}_{GPB}}-f{{\left\{ s,g \right\}}_{GPB}}$, can this equation get a solution for all $f$, it can be considered.  
With the properties given above, we can define a general form of Poisson manifolds.
\begin{definition}
The generalized structural Poisson bracket on the smooth manifold $M$ is an operation on the smooth real function space $ {{C}^{\infty }}\left( M,\mathbb{R} \right)$, for $f,g\in  {{C}^{\infty }}\left( M,\mathbb{R} \right)$, it exists $$\left\{ f,g\right\}=\left\{ f,g \right\}_{GPB}+G\left(s, f,g \right)\in  {{C}^{\infty }}\left( M,\mathbb{R} \right),$$ the operation satisfies the condition $(1)-(4)$, then $\left( M,S,\left\{ , \right\}\right)$ is said to be a generalized Poisson manifolds.
\end{definition}

\begin{definition}\label{d6}
Let $P$ be a Poisson manifold, geometric potential function $s:P\to \mathbb{R}$ is a smooth function. The potential vector field ${{X}_{s}}$ associated with $s$ is a smooth vector field on $P$, it satisfies $${{X}_{s}}\left( f \right)={{\left\{ f,s \right\}}_{GPB}}=-\hat{S}f,$$ for all smooth functions $f:P\to \mathbb{R}$, where $\hat{S}={{b}_{j}}{{\partial }_{j}}=-{{X}_{s}}$ is called structural operator,  then ${{X}_{s}}$ is the geometric potential vector field on $P$.
\end{definition}
It is clear to see that ${{X}_{s}}={{c}_{j}}{{\partial }_{j}}=-\hat{S}$, where $${{b} _{j}}={{J}_{ij}}{{A}_{i}}=-{{J}_{ji}}{{A}_{i}}=-{{c}_{j}},$$ and ${{c} _{j}}={{J}_{ji}}{{A}_{i}}$. It can be obtained from the structural operator  $G\left( s,s,f \right)=-G\left( s,f,s \right)=s\hat{S}f$.

\section{Generalized covariant Hamilton system}
In this section, we put the generalized structural Poisson bracket into practise, actually, it works for the another Hamilton system which is called the generalized covariant Hamilton system (GCHS). 

More precisely, the generalized covariant Hamilton system with respect to the function $f$ is expressed as $$\frac{\mathcal{D}f}{dt}=\left\{ f,H \right\}={{\left\{ f,H \right\}}_{GPB}}+G\left( s,f,H \right),$$ which naturally contains two parts as follows:
thorough generalized Hamiltonian system (TGHS)\newline~~~~~ $\frac{df}{dt}= {\dot {f}}={{\left\{ f,H \right\}}_{GPB}}-H{{\left\{ s ,f \right\}}_{GPB}}$,\newline
The S-dynamics:\newline $\frac{ds }{dt}=w={{\left\{ s ,H \right\}}_{GPB}}=\left\{ 1,H \right\}$,\newline
where $\frac{\mathcal{D}}{dt}=d/dt+w$ is the covariant time derivative.

The generalized covariant Hamilton system can also be expressed in the form of a vector field as
\begin{align}
  & \mathcal{D}f/dt=\left\{ f,H \right\}={{X}_{H}}\left( f \right)+H{{X}_{s}}\left( f \right)-f{{X}_{s}}\left( H \right) \notag\\
 & \begin{matrix}
   {} & {} & ={{X}_{H}}\left( f \right)+f\hat{S}H-H\hat{S}f,  \\
\end{matrix} \notag
\end{align}
The geometric bracket $$G\left( s,f,H \right)=H{{X}_{s}}\left( f \right)-f{{X}_{s}}\left( H \right)=f\hat{S}H-H\hat{S}f,$$ is completely induced by the potential vector field ${{X}_{s}}$. Obviously, the generalized covariant Hamiltonian system is divided into three major fields, namely, the conservative force field and two external force fields interacting with the manifold.

The generalized covariant Hamilton system in terms of the coordinates has generalized the generalized Hamilton system to the more general form
\begin{equation}\label{eq23}
  \frac{\mathcal{D}x}{dt}={\dot{x}}+wx,~~~x\in {{\mathbb{R}}^{n}},
\end{equation}
and it greatly opened up a new insight to study the Hamiltonian system with extra structure. The generalized covariant Hamilton system about coordinate is expressed as
$$\left\{ \begin{matrix}
   \mathcal{D}x/dt=J\left( x \right)\nabla H\left( x \right)+H\left( x \right)J\left( x \right)\nabla s\left( x \right)+\rho \left( x \right)  \\
   \dot{x}=J\left( x \right)DH\left( x \right)=J\nabla H+HJ\nabla s  \\
   w\left( x \right)=\dot{s}=\hat{S}H=\left\{ 1,H \right\}={{A}^{T}}J\nabla H  \\
\end{matrix} \right.$$where $\rho \left( x \right)=xw\left( x \right)$, obviously, the generalized covariant Hamilton system is composed of three parts, namely, the generalized Hamilton system $J\left( x \right)\nabla H\left( x \right)$ and the sum of two parts $HJ\nabla s\left( x \right)$ and $\rho \left( x \right)$ that interact with the manifold space or external environment. The formulations of the component are shown as follows,
$$\frac{\mathcal{D}{{x}_{k}}}{dt}=\left\{ {{x}_{k}},H\right\}={\dot{x}_{k}}+{{x}_{k}}w={{J}_{kj}}{{\partial}_{j}}H+H{{c}_{k}}+{{x}_{k}}w,$$
 $${\dot{x}_{k}}=\frac{d{{x}_{k}}}{dt}={{J}_{kj}}{{D}_{j}}H
 ={{J}_{kj}}{{\partial}_{j}}H+H{{c}_{k}},$$
\[w=\frac{ds}{dt}={{A}_{k}}{\dot{x}_{k}}={{b}_{j}}{{\partial }_{j}}H.\]Note that the S-dynamics describes the rotations of the manifold in an angular frequency.
As a result, there has equation of time evolution
$df/dt=\dot{f}={{\dot{x}}_{k}}{{\partial }_{k}}f$ for $f\in  {{C}^{\infty }}\left( M,\mathbb{R} \right)$, namely, the time derivative
$d/dt={{\dot{x}}_{j}}{{\partial }_{j}}$ is given, similarly, it gets covariant time derivative $\mathcal{D}/dt={{\dot{x}}_{k}}{{D}_{k}}=d/dt+w$. It implies that the ${\dot{x}_{k}}={{J}_{kj}}{{D}_{j}}H$ is the foundation for the time evolution in terms of a function. In some aspects, the $x,s$ are the primary variables on the manifolds to give the different dynamical explanations.

Let's recall the conservative force, it says that a particle of mass $m$ moves under the influence of a conservative force derived from the gradient $\nabla$ of a scalar potential, the force is given by  $F =-\nabla V(x)$, where $V$ is potential energy, and  $F_{j} =-\partial_{j} V(x)$. As the generalized gradient operator given, then it leads to the consideration of application to the potential energy, then the gradient is replaced by the generalized gradient, and $F =-\nabla V(x)-VA$, or $F_{j} =-\partial_{j} V(x)-VA_{j}$, according to $F=\frac{d}{dt}p$, it has a general form given by
\begin{equation}\label{a2}
  \frac{d}{dt}p_{j}=F_{j} =-\partial_{j} V(x)-VA_{j}=-D_{j}V(x).
\end{equation}
Note that the gradient system can be replaced by the generalized gradient for a extensive generalization of corresponding system.
Using the generalized structural Poisson bracket to the functions $f,V$, it gets 
$$\left\{ f,V \right\}={{\left\{ f,V \right\}}_{GPB}}+G\left( s,f,V \right)={{\left\{ f,V \right\}}_{GPB}}-V{{\left\{ s,f \right\}}_{GPB}}+f{{\left\{ s,V \right\}}_{GPB}},$$where
the geometric bracket in terms of functions $f,V$ is given by \[G\left( s,f,V \right)=f{{\left\{ s,V \right\}}_{GPB}}-V{{\left\{ s,f \right\}}_{GPB}}.\]
If let $f=p_{j}$ be taken into above equation, then
\[\left\{ {{p}_{j}},V \right\}={{\left\{ {{p}_{j}},V \right\}}_{GPB}}-V{{\left\{ s,{{p}_{j}} \right\}}_{GPB}}+{{p}_{j}}{{\left\{ s,V \right\}}_{GPB}}.\]
In comparison with the \eqref{a2}, it implies that 
 ${{A}_{j}}={{\left\{ s,{{p}_{j}} \right\}}_{GPB}}$.  
Conclusively, it's convinced that the classical Poisson bracket and classical Hamilton system, the generalized Poisson bracket and generalized Hamiltonian system, the generalized structural Poisson bracket and generalized covariant Hamiltonian are in pairs working for the physical system.

Within the framework of the generalized covariant Hamilton system, there are two second order equations of the generalized covariant Hamilton system with respect to time $\frac{{{\mathcal{D}}^{2}}{}}{d{{t}^{2}}}$ and coordinates ${{D}_{k}}\frac{\mathcal{D}{{x}_{k}}}{dt}$, it will be discussed separately,

\begin{theorem}\label{t2}
There exists an equality
 $D\cdot \frac{\mathcal{D}x}{dt}=\nabla\cdot v+\left( x\cdot D+2\right)w$.

\begin{proof}Based on the equation (\ref{eq23}), let $D_{l}$ act on the component expression of generalized covariant Hamilton system in terms of the coordinates and obtain
\begin{align}\label{a1}
 {{D}_{l}}\frac{\mathcal{D}{{x}_{k}}}{dt}
  &={{D}_{l}}{\dot{x_{k}}}+{{D}_{l}}\left( w{{x}_{k}} \right)
={{D}_{l}}{\dot{x}_{k}}+{{x}_{k}}{{D}_{l}}w+w{{\delta }_{lk}}   \\
 & ={{\partial }_{l}}{\dot{x}_{k}}+{{A}_{l}}{\dot{x}_{k}}+{{x}_{k}}{{\partial }_{l}}w+{{x}_{k}}{{A}_{l}}w+w{{\delta }_{lk}}. \notag
\end{align}
If $l=k$ holds for the \eqref{a1}, then the consequence will be accordingly obtained
$${{D}_{k}}\frac{\mathcal{D}{{x}_{k}}}{dt}=\nabla \cdot v+\nabla s \cdot v+\left( x\cdot D+1 \right)w=\nabla \cdot v+\left( x\cdot D+2 \right)w,$$
where the divergence of vector $v$ is $\nabla \cdot v={{\partial }_{k}}{\dot{x}_{k}}$ and $w=\nabla s \cdot v=A \cdot v$, where ${{v}_{k}}={\dot{x}_{k}}$ has been used.
\end{proof}
\end{theorem}
As a consequence, if $\frac{\mathcal{D}x}{dt}=0$, then $\nabla \cdot v+ x\cdot Dw=-2w$.

Let's consider the formula ${{D}_{l}}w$, by computing, it gets
\[{{D}_{l}}w={{D}_{l}}\left( {{A}_{k}}{{\dot{x}_{k}}} \right)={{\dot{x}_{k}}}{{B}_{kl}}+{{A}_{k}}{{\partial }_{l}}{{\dot{x}_{k}}},\]
where ${{B}_{kl}}={{\partial }_{l}}{{A}_{k}}+{{A}_{k}}{{A}_{l}}={{D}_{l}}{{A}_{k}}$ can be realized as the curvature induced by the geometric potential function.

When we consider the second order form of the generalized covariant Hamilton system for the function $f\in {{C}^{\infty }}\left( M,\mathbb{R} \right)$, it naturally appears that
\[\frac{{{\mathcal{D}}^{2}}}{d{{t}^{2}}}f=\ddot{f}+2w\dot{f}+\beta f,\]holds for all functions $f$, where $\frac{df}{dt}=\dot{f},\frac{{{d}^{2}}f}{d{{t}^{2}}}=\ddot{f}$, and $\beta ={{w}^{2}}+\frac{d}{dt}w$. Apparently, this equation simultaneously contains two kinetic quantities $\dot{f},w=\dot{s}$. In fact, we can denote operator $\hat{\alpha }=2w\frac{d}{dt}+\frac{d}{dt}w$, then the second order covariant derivative of time is 
$\frac{{{\mathcal{D}}^{2}}}{d{{t}^{2}}}
=\frac{{{d}^{2}}}{d{{t}^{2}}}+{{w}^{2}}+\hat{\alpha }$,and $\hat{\alpha }$ satisfies an identity $\hat{\alpha }w^{-1/2}\equiv 0$.  
Obviously, when it comes to the coordinates, it gets
$${{a}_{i}}=\frac{{{\mathcal{D}}^{2}}{{x}_{i}}}{d{{t}^{2}}}
={\ddot{x}_{i}}+2w{\dot{x}_{i}}+{{x}_{i}}\beta,$$
that can be regarded as general geodesic equation, where
${\ddot{x}_{i}}=\frac{{{d}^{2}}}{d{{t}^{2}}}{{x}_{i}}$,
due to $w=\sum\limits_{j}{{{A}_{j}}{\dot{x}_{j}}}={{A}_{j}}{\dot{x}_{j}}$, then
${{a}_{i}}={\ddot{x}_{i}}+2A_{j}{\dot{x}_{j}}{\dot{x}_{i}}+{{x}_{i}}\beta$. As it can be seen, the term $2w{\dot{x}_{i}}+{{x}_{i}}\beta$ is associated with the S-dynamics, this feature indicates the importance of the S-dynamics.
Therefore, if ${{a}_{i}}=0$, then it yields
$${\ddot{x}_{i}}+2A_{j}{\dot{x}_{j}}{\dot{x}_{i}}=-\beta{{x}_{i}}. $$ Furthermore, if
$\beta=0$ is assumed, then it gets a special solution of the S-dynamics $w=\frac{1}{t+{{C}_{1}}}$, where ${C}_{1}$ is a integral constant.
${\ddot{x}_{i}}+\frac{2}{t+{{C}_{1}}}{\dot{x}_{i}}=0$

\subsection{Example}
In this section, we give a simple example to show the generalized covariant Hamilton system in terms of the coordinates. 
Consider differential equations of motion
\begin{align}
  & {\dot{x}_{1}}={{x}_{2}}, \notag\\
 & {\dot{x}_{2}}=-{{x}_{1}},\notag\\
 & {\dot{x}_{3}}=-4{{x}_{1}}{{x}_{2}}. \notag
\end{align}
The Hamiltonian function of the system is 
$H=\frac{1}{2}\left( x_{1}^{2}+x_{2}^{2} \right)$ and then it has ${{\left( {{\partial }_{1}}H,{{\partial }_{2}}H,{{\partial }_{3}}H \right)}^{T}}={{\left( {{x}_{1}},{{x}_{2}},0 \right)}^{T}}$, the structural matrix is
\[\left( {{J}_{ij}} \right)=\left( \begin{matrix}
   0 & 1 & 2{{x}_{2}}  \\
   -1 & 0 & 2{{x}_{1}}  \\
   -2{{x}_{2}} & -2{{x}_{1}} & 0  \\
\end{matrix} \right).\]
Its generalized covariant Hamilton system is divided into two parts.  First, consider the thorough generalized Hamilton system, 
$${{\dot{x}}_{1}}={{x}_{2}}-H{{b}_{1}},$$
$${{\dot{x}}_{2}}=-{{x}_{1}}-H{{b}_{2}},$$
$${{\dot{x}}_{3}}=-4{{x}_{1}}{{x}_{2}}-H{{b}_{3}},$$
where $$\left(\begin{matrix}
   {{c}_{1}}  \\
   {{c}_{2}}  \\
   {{c}_{3}}  \\
\end{matrix} \right)=-\left(\begin{matrix}
   {{b}_{1}}  \\
   {{b}_{2}}  \\
   {{b}_{3}}  \\
\end{matrix} \right)=\left( \begin{matrix}
   {{A}_{2}}+2{{x}_{2}}{{A}_{3}}  \\
   -{{A}_{1}}+2{{x}_{1}}{{A}_{3}}  \\
   -2{{x}_{2}}{{A}_{1}}-2{{x}_{1}}{{A}_{2}}  \\
\end{matrix} \right).$$
The corresponding S-dynamic system is calculated as 
$$w={{x}_{2}}{{A}_{1}}-{{x}_{1}}{{A}_{2}}-4{{x}_{1}}{{x}_{2}}{{A}_{3}}.$$
For specific calculation, the geometric scalar potential function is assumed to be
$s=\ln \left( x_{1}^{2}+x_{2}^{2}+x_{3}^{2} \right)=\ln {{l}^{2}}$, where 
${{l}^{2}}=x_{1}^{2}+x_{2}^{2}+x_{3}^{2}$, then the structural derivative is
$${{\left( {{A}_{1}},{{A}_{2}},{{A}_{3}} \right)}^{T}}=2{{\left( {{x}_{1}},{{x}_{2}},{{x}_{3}} \right)}^{T}}/{{l}^{2}},$$ then 
$${{c}_{1}}=2\left( {{x}_{2}}+2{{x}_{2}}{{x}_{3}} \right)/{{l}^{2}},{{c}_{2}}=\left( -2{{x}_{1}}+4{{x}_{1}}{{x}_{3}} \right)/{{l}^{2}},{{c}_{3}}=-8{{x}_{1}}{{x}_{2}}/{{l}^{2}}.$$
The integrated S-dynamics is 
$w=-8{{x}_{1}}{{x}_{2}}{{x}_{3}}/{{l}^{2}}$.

\section{Conclusions}
As the generalized structural Poisson bracket proposed, it reminds us of the Lie bracket, there should have a larger range that Lie bracket can not reach, for this reason, based on non-Abelian Lie bracket, the generalized geometric Lie bracket on vector space is proposed for a further explanation, that is the Lie bracket equipped with the geometric bracket on the vector space.  Then we briefly discuss the  gradient system related to the conservative force, and some results in second order of the generalized covariant Hamilton system follow.

\end{document}